# Proof of the Goldbach conjecture in a more stringent form.

Jürgen Schwarz

Institut für Physikalische Chemie, Arnold-Sommerfeld-Str. 4
38678 Clausthal-Zellerfeld, Germany

e-mail: sekretariat@pc.tu-clausthal.de

## Abstract

With an artificial (*p*', *n*')-system it has been proved that even numbers > $p_x^2$ are the sum of two $p > p_x$.

To prove the conjecture of Goldbach, that all even numbers > 2 are a sum of two prime numbers, we use an artificial ($p'$, $n'$) “fantom” system, in which the $p'$ are not necessarily prime ($p$) and the $n'$ not necessarily nonprime ($n$).

A fantom system $F(p_x)$ is related to a $p_x$, where $x$ is the serial index of the $p$. In the length $L(p_x) = \prod_{i=1}^{x} p_i$ of the natural numbers all $p$ from 2 to $p_x$ and the products containing them are cancelled. The cancelled numbers are $n'$, the remaining numbers are $p'$. The distribution of the $p'$ in the length $L(p_x)$ is symmetric.

The first three of the infinite series of fantom systems are

$F(2)$: 1 $L(2) = 2$;

$F(3)$: 1, 5 $L(3) = 6$;

$F(5)$: 1, 7, 11, 13, 17, 19, 23, 29 $L(5) = 30$.

An $F(p_x)$ can be constructed by lining up $p_x$ systems $F(p_{x\text{-}1})$ to a presystem $PF(p_x)$ and then canceling the products of $p_x$ with the $p'$ of $F(p_{x\text{-}1})$:

```
F(2) → F(3):  1, 1, 1
              1, 3̸, 5

F(3) → F(5):  1, 5, 1,  5,  1,  5,  1,  5,  1,  5
              1, 5̸, 7, 11, 13, 17, 19, 23, 2̸5̸, 29.
```

Considering this transition from an $F(p_{x\text{-}1})$ to $F(p_x)$ leads to $A(p_x) = \prod_{i=1}^{x} [p_i - 1]$ for the amount of the $p'$ in an $F(p_x)$.

If we multiply the series of $p'$ of a fantom system $F(p_x)$ with one of these $p'$, we get a new series of numbers of the form $a \cdot L(p_x) + p'$, which are distributed symmetrically over the appropriate length, and $a$ is an integer including 0, and $p'$ is one of the series numbers. We give an example just for the change of the series of the $p'$ for $F(5)$:

```
(1, 7, 11, 13, 17, 19, 23, 29) · 7
→ (7, 19, 17, 1, 29, 13, 11, 23)
```

There are for this system 8 or generally $A(p_x)$ possible series.

The fantom systems and the resulting systems of sums are highly regular. The regularities in any deliberately chosen range of a system are valid also in the range. where the $p' = p$ and the $n' = n$. This is the basis of the following proof.

We define a set of sums of a fantom system $S(p_x)$ as all sums of any two $p'$ of $F(p_x)$, where sums with different summands are counted twice. So we have $[A(p_x)]^2$ sums in $S(p_x)$, which are distributed in 2 $L(p_x)$. The distribution is symmetrical. We define a reduced set of sums $RS(p_x)$, in which all sums > $L(p_x)$ are replaced by the respective sums minus $L(p_x)$. This set is distributed in $L(p_x)$. If we line up two of these systems, where the second is shifted by $L(p_x)$ against the first, we see, that the sums are distributed symmetrically around $L(p_x)$.

```
S(2): |. 2 . .|                          RS(2):  |. 2|
S(3): |. 2 . . . 6 . . . 10 . .|         RS(3):  |. 2 . 4 . 6|
                 6                                          6
```

An $RS(p_x)$ can be constructed by arranging $p_x^2$ systems $RS(p_{x-1})$ in a quadratic form to a presystem $PRS(p_x)$ and then canceling the sums, which contain as a summand a product of $p_x$ with the $p'$ of $F(p_{x-1})$.

| | 2 | 4 | 6 |
|---|---|---|---|
| 1 | 1 | ~~3~~ | 5 |
| 3 | ~~5~~ | ~~1~~ | ~~3~~ |
| 5 | ~~3~~ | 5 | 1 |
| | 1 | 1 | 2 |

*PRS*(3) → *RS*(3)

| | 2 | 4 | 6 | 8 | 10 | 12 | 14 | 16 | 18 | 20 | 22 | 24 | 26 | 28 | 30 |
|---|---|---|---|---|---|---|---|---|---|---|---|---|---|---|---|
| 1 | 1 | | ~~5~~ | 7 | | 11 | 13 | | 17 | 19 | | 23 | ~~25~~ | | 29 |
| 3 | | | | | | | | | | | | | | | |
| 5 | | ~~29~~ | ~~1~~ | | ~~5~~ | ~~7~~ | | ~~11~~ | ~~13~~ | | ~~17~~ | ~~19~~ | | ~~23~~ | ~~25~~ |
| 7 | ~~25~~ | | 29 | 1 | | ~~5~~ | 7 | | 11 | 13 | | 17 | 19 | | 23 |
| 9 | | | | | | | | | | | | | | | |
| 11 | | 23 | ~~25~~ | | 29 | 1 | | ~~5~~ | 7 | | 11 | 13 | | 17 | 19 |
| 13 | 19 | | 23 | ~~25~~ | | 29 | 1 | | ~~5~~ | 7 | | 11 | 13 | | 17 |
| 15 | | | | | | | | | | | | | | | |
| 17 | | 17 | 19 | | 23 | ~~25~~ | | 29 | 1 | | ~~5~~ | 7 | | 11 | 13 |
| 19 | 13 | | 17 | 19 | | 23 | ~~25~~ | | 29 | 1 | | ~~5~~ | 7 | | 11 |
| 21 | | | | | | | | | | | | | | | |
| 23 | | 11 | 13 | | 17 | 19 | | 23 | ~~25~~ | | 29 | 1 | | ~~5~~ | 7 |
| 25 | ~~7~~ | | ~~11~~ | ~~13~~ | | ~~17~~ | ~~19~~ | | ~~23~~ | ~~25~~ | | ~~29~~ | ~~1~~ | | ~~5~~ |
| 27 | | | | | | | | | | | | | | | |
| 29 | | ~~5~~ | 7 | | 11 | 13 | | 17 | 19 | | 23 | ~~25~~ | | 29 | 1 |
| | 3 | 3 | 6 | 3 | 4 | 6 | 3 | 3 | 6 | 4 | 3 | 6 | 3 | 3 | 8 |

*PRS*(5) → *RS*(5)

In this diagrams any field represents a sum of the number in the field and the number in the same height to left outside. Above the diagrams are the sums, the even numbers, below the amounts of the sums.

Some general facts should be seen. The amount of sums, which give a certain even number equals to the amount of different summands, which are involved, if sums of different summands are counted twice.

The symmetry of the amounts of sums belonging to certain even numbers, which we have seen above for the $RS(p_x)$, leads to an internal symmetry in the $RS(p_x)$ in the range from 2 to $L(p_x) - 2$: The amounts of sums, which belong to two even numbers, which add up to $L(p_x)$ is equal. This is the case for the even numbers 2 and 4 in $RS(3)$ or for 2 and 28, 4 and 26, and so on in $RS(5)$. The even number $L(p_x)$ is always built by an extraordinary high amount of sums, since all $p'$ of $F(p_x)$ are involved as summands.

We will consider the diagrams in more detail. In $PRS(3)$ are canceled all sums containing $3 \cdot 1 = 3$ by a horizontal line, associated with a diagonal line, which is broken in two parts. In $PRS(3)$ the summands 1, $2 + 1 = 3$, $2 \cdot 2 + 1 = 5$ are involved in all sums 2, $2 + 2 = 4$, $2 \cdot 2 + 2 = 6$. One should realize, that 2 results from $3 + 5 = 8$; $8 - 6 = 2$ and 4 from $5 + 5 = 10$; $10 - 6 = 4$. One should see, that the amount of canceling in 2, 4 and 6, which is 2, 2 and 1, is a consequence of the fact, that in $F(2)$ is one $p'$, the number 1, which leads to the one canceling number 3, which works once in the sum with itself in 6, or works twice in sums with different summands in 2 and 4.

In $PRS(5)$ are cancelled all sums containing 5 and 25 by two horizontal lines, each associated with a diagonal line, which is broken in two parts. The canceling summands 5 and 25 follow from $(1, 5) \cdot 5 \rightarrow (5, 4 \cdot 6 + 1)$. The sums 2, $6 + 2 = 8$, $2 \cdot 6 + 2 = 14$, $3 \cdot 6 + 2 = 20$, $4 \cdot 6 + 2 = 26$ are built from the summand 1 and its equivalents $6 + 1 = 7$, $2 \cdot 6 + 1 = 13$, $3 \cdot 6 + 1 = 19$, $4 \cdot 6 + 1 = 25$. In any of these sums, the amount of which is 5 in $PRS(5)$, the amount is reduced by the two sums with 25, except in 20 where the horizontal and the diagonal line intersect and a sum of two cancelling summands is built. One should see, that the reductions of 2 or 1 are the consequence of the fact, that in the single sum 2 in $RS(3)$ only one summand is involved. The sums 4, $6 + 4 = 10$, $2 \cdot 6 + 4 = 16$, $3 \cdot 6 + 4 = 22$, $4 \cdot 6 + 4 = 28$ are built from the summand 5 and its equivalents $6 + 5 = 11$, $2 \cdot 6 + 5 = 17$, $3 \cdot 6 + 5 = 23$, $4 \cdot 6 + 5 = 29$. In any of these sums, the amount of which is 5 in $PRS(5)$, the amount is reduced by two sums with 5, except in 10, where a sum of two cancelling summands is built. The reduction of 2 or 1 are the consequence of the fact, that in the single sum 4 in $RS(3)$ only one summand is involved. The sums 6, $2 \cdot 6$, $3 \cdot 6$, $4 \cdot 6$, $5 \cdot 6$ are built from the summands 1 and 5 and their equivalents. In any of these sums, the amount of which is 10 in $PRS(5)$, the amount is reduced by the two sums with 25 and the two sums with 5, except 30, where we have not 4 summands involved in the canceled sums, but only 2 because of intersection. The reductions of 4 or less are the consequence of the fact, that in the double sum 6 in $RS(3)$ two summands are involved. For $PRS(7) \rightarrow RS(7)$ the description and argumentation is principally analog. The canceling summands follow from $(1, 7, 11, 13, 17, 19, 23, 29) \cdot 7$. The sums 2, 30

+ 2, 2 · 30 + 2, … 6 · 30 + 2 are built from the summands 1, 19, 13 and their equivalents. Three of these 21 summands are canceling summands.

We can generalize for the transition of $PRS(p_x) \rightarrow RS(p_x)$: In any length $L(p_{x-1})$ of the even numbers of $PRS(p_x)$ are canceled at most 2 $RS(p_{x-1})$. Since in any length $L(p_{x-1})$ of $PRS(p_x)$ are $p_x$ $RS(p_{x-1})$, there remain at least $[p_x - 2]$ $RS(p_{x-1})$.

A balance of the amounts of sums shows: The amount of sums in $PRS(p_x)$ is $[A(p_{x-1})]^2 \cdot p_x^2$. The amount of sums in $RS(p_x)$ is

$$[A(p_x)]^2 = [A(p_{x-1}) \cdot [p_x - 1]]^2 = [A(p_{x-1})]^2 \cdot p_x^2 - 2\, p_x\, [A(p_{x-1})]^2 + [A(p_{x-1})]^2$$

The change in the amount of sums from $PRS(p_x)$ to $RS(p_x)$ is in any of the $p_x$ lengths $L(p_{x-1})$ of the even numbers of the system $-2\,[A(p_{x-1})]^2 + \varepsilon$, where $\varepsilon$ is a relatively small number of sums. All $\varepsilon$ values, which are distributed over the system, add up to one $[A(p_{x-1})]^2$. They stem from the sums between canceling summands, the amount of which is $A(p_{x-1})$.

In $RS(3)$ the $\varepsilon$ is just the one value 6 = 2 · 2 + 2, which is equivalent to $RS(2)$. In $RS(5)$, the $\varepsilon$ are 10 = 6 + 4, 20 = 3 · 6 + 2 and two times 30 = 4 · 6 + 6, which is equivalent to $RS(3)$.

So we have proved by induction, that any even number in an $RS(p_x)$ can be built by an amount of at least $\prod_{i=2}^{x} [p_i - 2]$ sums, where sums of different summands are counted twice.

In a fantom system $F(p_x)$ are the $p' = p$ and the $n' = n$ in the range $> p_x$ and $< [p_{x+1}]^2$. In the lower range $1 = p'$ and the $p_i$ from 2 to $p_x$ are $n'$. In the higher range the products between $p_i > p_x$ are $p'$. The smallest product of this kind is $[p_{x+1}]^2$.

It is appropriate to consider the range from $p_x^2$ to $[p_{x+1}]^2$ in the conditions of $F(p_x)$ and then the range from $[p_{x+1}]^2$ to $[p_{x+2}]^2$ in the conditions of $F(p_{x+1})$ and so on.

We consider the situation in the framework of the diagrams of $RS(p_x)$, from which the two with $x$ = 2 and 3 are shown. The construction of the diagrams is clearly defined.

The above mentioned range is principally on the right side of the diagonal which divides sums and reduced sums.

The minimum amount of sums for any even number is $\prod_{i=2}^{x} [p_i - 2]$, the amount of places is ½ $L(p_x) = \prod_{i=2}^{x} p_i$ . So we have the density as the quotient $\prod_{i=2}^{x} [1 - 2 p_i^{-1}]$. The length $p_x^2$ occupies ½ $[p_x^2 + 1]$ places, which multiplied with the density lead to an expression of a mean amount of sums in the length $p_x^2$. This amount increases with $x$.

We have to know however from an independent consideration the maximum of possible fluctuations in the amount of sums for an even number in a length $L$, which is chosen deliberately within the length $L(p_x)$. The superposition of $x$ different combs, the teeth of which are separated by $p_i$ with $i$ from 1 to $x$, shows the situation of canceled numbers in $F(p_x)$. The variance in the number of teeth within $L$ is 1 for each comb, except for the case, that $L/p_i$ is an integer, where the variance is 0. So we have the maximum variance $x$ as well for canceled as for not canceled numbers. This is independent from the chosen length $L$.

If we move $L$ step by step over the system of superposed combs, at the entrance of a tooth, single or superposed with up to $x – 1$ teeth, a tooth of the related comb leaves $L$ at the other end at the same time, if the variance is 0. Otherwise it leaves $L$, before the next tooth of this kind enters $L$. If an $L > x$ is filled completely with teeth, then in the extreme opposite case there are at most $x$ free places. $x$ is the number of teeth, which can be inserted deliberately, whether superposed on one place or lined up or somehow.

The situation of canceled sums belonging to an even number is described by the situation of superposition of two such sets of combs, which are all shifted against each other deliberately without matching of combs of the same kind, because two summands are involved. So we have the maximum variance $2x$.

The final formula for the minimum amount of sums, which build the even number $p_x^2 + 1$ is $C = \frac{1}{2}\prod_{i=2}^{x}\left[1-2p_i^{-1}\right] \cdot \frac{1}{2}\left[p_x^2+1\right] - 1 - 2x + 2$.

The first $\frac{1}{2}$ stands, because in $RS(p_x)$ the sums of different summands are counted twice. The –1 stands for the possible sum with the summand 1. The $–2x$ is the subtraction of the maximum variance, which is reduced by 2, because we consider even numbers. $C$ contains not sums with summands $p_i$ from 3 to $p_x$. $C$ increases, if we go through the even numbers from $p_x^2 + 1$ to $[p_{x+1}]^2 – 1$. Numerical investigation shows that $C$ becomes > 1 at $53^2 + 1 = 2810$. It diverges for $x \to \infty$.

We have proved for numbers > 2810 the Goldbach conjecture in the more stringent form, that all even numbers $> p_x^2$ are the sum of two $p_i > p_x$. We suppose, that this more stringent form is valid for all even numbers > 4.